\documentclass[11pt,reqno]{amsart}
\usepackage{amsmath,amssymb,graphicx,amscd,hyperref, amsfonts,psfrag,fancybox,qsymbols,indentfirst}
\usepackage{mathrsfs}
\usepackage{dsfont}
\usepackage[latin1]{inputenc}
\renewcommand{\S}{\mathscr{S}}  

\newcommand{\NN}{{\mathbb{N}}}

 \DeclareMathOperator{\E}{\mathbb{E}}      
 \DeclareMathOperator{\Var}{\mathrm{Var}}      
 \DeclareMathOperator{\cov}{cov}

\DeclareMathOperator{\W}{Wg}

 \newcommand{\ii}{{\mathrm{i}}}
  
\newcommand{\law}{\overset{\mbox{\rm \scriptsize law}}{=}}

\newcommand{\convlaw}{\overset{\mbox{\rm \scriptsize law}}{\longrightarrow}}
\newcommand{\convprob}{\overset{\mbox{\rm \scriptsize P}}{\longrightarrow}}
\newcommand{\convfidi}{\overset{\mbox{\rm \scriptsize fidi}}{\longrightarrow}}

\newtheorem{thm}{Theorem}[section]

\newtheorem{lem}[thm]{Lemma}
\newtheorem{prop}[thm]{Proposition}
\newtheorem{rem}[thm]{Remark}
\theoremstyle{definition}

\theoremstyle{remark}

\def \be{\begin{eqnarray*}}
\def \ee{\end{eqnarray*}}
\def \ben{\begin{eqnarray}}
\def \een{\end{eqnarray}}

\def\sn{^{(n)}}

\numberwithin{equation}{section}

\begin{document}
\title[Random truncation  and  bridges ]
{Random truncations of  Haar distributed matrices and  bridges}
\date{\today}

\author{C. Donati-Martin}
\address{Universit\'e  Versailles-Saint Quentin, LMV UMR 8100,
 B\^atiment Fermat, 45 avenue des Etats-Unis,
F-78035 Versailles Cedex}
\email{catherine.donati-martin@uvsq.fr}
\urladdr{http://lmv.math.cnrs.fr/annuaire/catherine-donati/}

\author{A. Rouault}
 \address{Universit\'e  Versailles-Saint Quentin, LMV UMR 8100,
 B\^atiment Fermat, 45 avenue des Etats-Unis,
F-78035 Versailles Cedex}
 \email{alain.rouault@uvsq.fr}
\urladdr{http://rouault.perso.math.cnrs.fr/}

\subjclass[2010]{15B52,  60F17, 60J65} 
\keywords{Random Matrices, unitary ensemble, orthogonal ensemble, bivariate Brownian bridge, 
 invariance principle}

\maketitle

\begin{abstract} Let $U$ be a Haar distributed
 matrix in $\mathbb U(n)$ or $\mathbb O (n)$.  In a previous paper, we proved that after centering, the two-parameter process 
\[T^{(n)} (s,t) = \sum_{i \leq \lfloor ns \rfloor, j \leq \lfloor nt\rfloor} |U_{ij}|^2\]
converges in distribution to the bivariate tied-down Brownian bridge. In the present paper, we replace the deterministic truncation of $U$ by a random one, where each row (resp. column) is chosen with probability $s$ (resp. $t$) independently. We prove that the corresponding  two-parameter process, after centering and normalization by $n^{-1/2}$ converges to a Gaussian process. On the way we meet other interesting convergences.
\end{abstract}

\section{Introduction}
Let us consider  a   unitary matrix $U$ of size $n\times n$.  We fix two integers $p< n$ and $q < n$ and delete deterministically $n- p$ rows and $n- q$ columns.  Let us call $U^{p,q}$ the (rectangular) matrix so obtained. It is well known that if $U$ is Haar distributed in $\mathbb U(n)$, the random matrix $U^{p,q}\left(U^{p,q}\right)^*$ has a Jacobi matricial distribution and that if $(p/n, q/n) \rightarrow (s,t) \in(0,1)^2$, 
 its empirical spectral distribution converges to a limit ${\mathcal D}_{s,t}$ (see for instance \cite{Collins}), often called the generalized Kesten-McKay distribution. It is clear from the invariance of the Haar distribution on $\mathbb U (n)$ that we can keep the first $p$ rows and the first $q$ columns. 

In \cite{CDMAR} we studied the trace of  $U^{p,q}\left(U^{p,q}\right)^*$. It is the squared of the  Frobenius (or Euclidean) norm of $U^{p,q}$. Actually we set $p =\lfloor ns\rfloor, q= \lfloor nt\rfloor$ and considered the process indexed by $s,t \in [0,1]$. We proved that, after centering, but without any normalization, the process converges in distribution, as $n \rightarrow \infty$ to a bivariate tied-down Brownian bridge. Previously, Chapuy \cite{chapuy} proved a similar result for permutation matrices, with a $n^{-1/2}$ normalization. 

Besides, for purposes of random geometry analysis, in \cite{Farfourier} (see also \cite{Fararxiv}) B. Farrell 
deletes randomly and independently a proportion $1- s$ of rows and a proportion $1- t$ of columns. Let us call $\mathcal U^{s,t}$ the matrix so obtained. Farrell
  proved that (for fixed $s,t$) the empirical spectral distribution of $\mathcal U^{s,t}\left(\mathcal U^{s,t}\right)^*$  has an explicit  limiting distribution which is precisely 
 ${\mathcal D}_{s,t}$. 
Actually,  Farrell considered first  the (deterministic) DFT matrix  
\[F_{jk} = \frac{1}{\sqrt n} e^{-2\ii \pi (j-1)(k-1)/n} \ , j,k = 1, \dots, n\,,\]
and proved  that a Haar unitary matrix  has the same behaviour.

It is then tempting to study the same statistic as above, when using this random truncation. Of course, instead of the DFT, we can as well consider any (random or not)
 matrix whose all elements are of modulus $n^{-1/2}$, for instance a (normalized) complex Hadamard  matrix.

To define all random truncations simultaneously and get  a two-parameter process, we define a double array of $n^2$ auxiliary uniform variables, in order to generate Bernoulli variables. We will prove below that after centering, we need a normalization to get a Gaussian limiting process. 

 We use the 
 Skorokhod space
  $D([0,1]^2)$.  It consists of functions
from $[0, 1]^2$ to $\mathbb R$ which are at each point right continuous (with respect
to the natural partial order of $[0, 1]^2$) and admit limits in all "orthants". 
The space  $D([0,1]^2)$ is endowed with the topology of Skorokhod (see \cite{bickel1971convergence} for the definition).

For the sake of completeness, we treat also the 
 one-parameter process, i.e. truncation of the first column of the unitary matrix, and the case of permutation matrices.

In Farell \cite{Farfourier}  is also mentioned another way of truncation, used in \cite{Tul}, consisting in  drawing only one array of $n$ Bernoulli variables to determine the choice of rows and columns. We did not consider this model here, to make the paper shorter.

The rest of the paper is organized as follows. In Sec.  2 we introduce the basic notations with the  main statistics and the list of limiting processes.  In Sec. 3 we present the statements of convergence in distribution for all models. In particuler   Theorem \ref{twodim} says that at the first order, the difference between DFT, Haar unitary or Haar orthogonal random matrices is not seen by the statistics $\mathcal T\sn$ at the limit since only the randomness of the truncation contributes. Besides,  Proposition \ref{twodimprop} says that this difference is seen at the second order, by means of the process $\mathcal Z\sn$. Sec. 4 is devoted to the proofs. The most delicate part is the proof of Prop. \ref{twodimprop} (2).  We first prove the convergence of finite dimensional distributions after replacement of the indicator variables by Gaussian ones (Lindeberg's strategy). Then we  prove tightness of the process rescaled by $n^{-1/2}$ by application of  Davydov and Zitikis's criterion  \cite{DaZ}. 
 It could be noted that before that, we tried  several tightness critera for the unscaled  (two-parameter)  process, such as Bickel and Wichura \cite{bickel1971convergence} or Ivanoff \cite{Iv}, but they  failed\footnote{We hope to address a definite answer to the question of  weak convergence of  $\mathcal Z\sn$ in a forthcoming paper.}.  In Sec. 6, we gather all the estimates of moments; i.e  of polynomial  integrals on $\mathbb U (n)$ and $\mathbb O(n)$. 
Finally in Sec. 7, we give the proof of the Lindeberg's strategy.

\section{Notations}
We introduce the random processes that we will consider in this paper and the various limiting processes involved. 
\subsection{The main statistics}
Let $U \in \mathbb U(n)$  be the unitary group of size $n$ and  let $U_{i,j}$ 
be the generic element of $U$.  Let now
$R_i, i=1, \cdots, n$ and $C_j, j=1, \cdots, n$ be two independent  families  of independent random variables uniformly distributed  on $[0,1]$.

For the one-parameter model, we introduce two processses. 
\begin{eqnarray}
B_0\sn &=& \left(\sum_1^{\lfloor ns\rfloor}\left( |U_{i1}|^2  -1/n\right)\ , \ s\in [0,1]\right),\\
{\mathcal B}_0\sn  &=& \left(\sum_1^n |U_{i1}|^2 \left(1_{R_i \leq s} - s\right)  \  , \ s\in [0,1]\right).
\end{eqnarray}

For the two parameters model, we introduce the same framework.
Let
\begin{equation}
\label{pq}T\sn_{s,t} = \sum_{i \leq \lfloor ns\rfloor, j\leq \lfloor nt\rfloor}  |U_{ij}|^2 \ , \end{equation}
and 
\[ T\sn = \left(T\sn_{s,t}, s,t \in [0,1]\right).\]
 Let now

\ben\label{st}{\mathcal T}\sn_{s,t}:= \sum_{i,j = 1}^n |U_{ij}|^2{\mathbf 1}_{R_i \leq s}{\mathbf 1}_{C_j \leq t} \ .
\een
We have clearly, since $|U_{ij}|^2$  and $(R_i, C_j)$ are independent, 
\begin{equation}
\label{st0}
 {\mathcal T}\sn_{s,t} - \mathbb E  {\mathcal T}\sn_{s,t} =  \sum_{i,j = 1}^n |U_{ij}|^2\left[{\mathbf 1}_{R_i \leq s}{\mathbf 1}_{C_j \leq t} -st\right].\end{equation}
Introducing the centering of Bernoulli variables, 
 and using the fact that the matrix $U$ is unitary, we may write 
\begin{equation}\label{anova} {\mathcal T}\sn_{st} - \mathbb E  {\mathcal T}\sn_{st} = {\mathcal Z}_{st}\sn + n^{1/2} \mathcal W_{s,t}\sn
\end{equation}
where 
\begin{eqnarray}
\label{zw}
{\mathcal Z}_{s,t}\sn &:=& \sum_{i,j = 1}^n |U_{ij}|^2\left[{\mathbf 1}_{R_i \leq s}-s \right]\left[{\mathbf 1}_{C_j \leq t} -t\right]\\
\mathcal W_{s,t}\sn &:=&sF_{t}\sn +  t G_{s}\sn
\end{eqnarray}
with the empirical processes
\begin{equation}
\label{FG}
F_{t}\sn := n^{-1/2} \sum_{j=1}^n [{\mathbf 1}_{C_j \leq t} - t] \ ; \ 
 G_{s}\sn := n^{-1/2} \sum_{i=1}^n [{\mathbf 1}_{R_i \leq s} - s]\,.
\end{equation}
In the sequel, we denote, 
\begin{eqnarray*}\mathcal T\sn = \left(\mathcal T\sn_{s,t}, s,t \in [0,1]\right) &,& \mathcal Z\sn = \left(\mathcal Z\sn_{s,t}, s,t \in [0,1]\right),\\
F\sn = \left(F\sn_t, t \in [0,1]\right) &,& G\sn = \left(G\sn_t, t \in [0,1]\right),
\end{eqnarray*}
and
\[\mathcal W\sn = \left(\mathcal W\sn_{s,t}, s,t \in [0,1]\right).\]

\subsection{Gaussian processes and bridges}

The classical Brownian bridge denoted by $B_0$ is a centered Gaussian process with continuous paths  defined on $[0,1]$ of covariance  
\[\mathbb E \left(B_0(s)B_0(s')\right) =  s\wedge s' - ss'\,.\] 
The bivariate Brownian bridge denoted by $B_{0,0}$  is a centered Gaussian process with continuous paths  defined on $[0,1]^2$ of covariance  
\[\mathbb E \left(B_{0,0}(s,t)B_{0,0}(s',t'))\right) =  (s\wedge s')(t\wedge t') - ss'tt'.\] 
The tied-down bivariate Brownian bridge denoted by $W^{(\infty)}$ is a centered Gaussian process with continuous paths  defined on $[0,1]^2$ of covariance  
\[\mathbb E [W^{(\infty)}(s,t) W^{(\infty)}(s',t')] = (s\wedge s'- ss')(t\wedge t' - tt').\]
Let also  ${\mathcal W}^{(\infty)}$ be the centered Gaussian process  with continuous paths  defined on $[0,1]^2$ of covariance 
\[\mathbb E [{\mathcal W}^{(\infty)}(s,t) {\mathcal W}^{(\infty)}(s',t')] =ss'(t\wedge t') + (s\wedge s')tt' -2ss'tt'.\]
It can be defined also as
\[{\mathcal W}^{(\infty)}(s,t) = sB_0(t) + t B'_0(s)\]
where $B_0$ and $B'_0$  two independent one parameter Brownian bridges.

At last we will meet the process denoted by $B_0\otimes B_0$ which is   a centered  process  with continuous paths  defined on $[0,1]^2$ by 
\[B_0 \otimes B_0 (s,t) = B_0(s) B'_0(t)\]
where $B_0$ and $B'_0$ are two independent Brownian bridges. This process is not Gaussian, but it has the same covariance as $W^{(\infty)}$.
\section{Convergence in distribution}
We present unfied results in the cases of unitary and orthogonal groups. For this purpose we use the classical notation
\[\beta' = \frac{\beta}{2}= \begin{cases}1/2 \;\;\;\;&\mbox{in the orthogonal case} ,\\
1\;\;&\mbox{in the unitary case.}
\end{cases}\]
Let us begin with the one-parameter processes, where $\convlaw$ means convergence in distribution in $D([0,1))$.
\begin{thm}
\label{onedim}
We have : 
\begin{eqnarray}
\label{pont1a}
\sqrt{n}B_0\sn \convlaw \sqrt{\beta'^{-1}}B_0,\\
\label{pont1b}
\sqrt{n}{\mathcal B}_0\sn \convlaw \sqrt{1+ \beta'^{-1}} B_0.
\end{eqnarray}
\end{thm}
\eqref{pont1a} is well known since at least Silverstein \cite{Silver1} (in the case $\beta'=1$). Both results are a direct consequence of the fact that the vector  $(|U_{i1}|^2, i=1, \dots, n)$ follows the   Dirichlet $(\beta', \dots, \beta')$ distribution on the simplex (see a detailed proof in Section 4).

Let us continue with  the  two-parameters processes, where  now $\convlaw$ means convergence in distribution in 
$D([0,1]^2)$, and $\convfidi$ means convergence in  distribution of finite dimensional marginals. The different normalizations are explained in Remark \ref{remn}.

We begin with a 
recall of the convergence in the deterministic model of truncation which  was proved in Theorem 1.1  of \cite{CDMAR}.

\begin{thm}
\label{cdmar}
Under the Haar measure on $\mathbb U (n)$ or $\mathbb O(n)$, 
\[\left(T\sn -\mathbb E T\sn\right) \convlaw \sqrt{\beta'^{-1}} W^{(\infty)}\,.\]
\end{thm}

\medskip

In the sequel we will prove the following convergences for the second model of truncation.
\begin{thm}
\label{twodim}
\begin{enumerate}
\item Under the DFT model, or more generally when $|U_{ij}|^2 =1/n$ a.s. for every $i,j$,  
\[n^{-1/2}\left(\mathcal T\sn - \mathbb E  \mathcal T\sn\right) \convlaw {\mathcal W}^{(\infty)}\,.\]
\item Under the Haar measure on $\mathbb U (n)$ or $\mathbb O(n)$, 
\[n^{-1/2}\left(\mathcal T\sn - \mathbb E  \mathcal T\sn\right) \convlaw {\mathcal W}^{(\infty)}\,.\]
\end{enumerate}
\end{thm}

The proof uses the decomposition
(\ref{anova}):
\[\mathcal T\sn - \mathbb E  \mathcal T\sn = n^{1/2} \mathcal W\sn + \mathcal Z\sn\,.\]  The following trivial lemma (consequence of the convergence of empirical processes $F\sn$ and $G\sn$) rules the behavior of $\mathcal W\sn$ and the following Proposition  \ref{twodimprop}  provides a convergence in distribution of $\mathcal Z\sn$  hence a convergence in probability to $0$ for $n^{-1/2}\mathcal Z \sn$.
\begin{lem}
\label{triv}
\begin{equation}
\label{convw}
\mathcal W\sn \convlaw \mathcal W^{\infty}\,.
\end{equation}
\end{lem}

\begin{prop}
\label{twodimprop}
\begin{enumerate}
\item Under the DFT model, or more generally if $|U_{ij}|^2 = 1/n$ a.s. for every $i,j$, 
\[\mathcal Z\sn \convlaw B_0\otimes B_0\,.\]
\item
Under the Haar model, 
\begin{equation} \label{maincv}
\mathcal Z\sn \convfidi  \sqrt{\beta'^{-1}}W^{(\infty)} + B_0\otimes B_0\,,
\end{equation}
 where  $B_0\otimes B_0$ and $W^{(\infty)}$  are independent, and
\ben
\label{zero}
n^{-1/2} \mathcal Z\sn \convprob 0\,.
\een
\end{enumerate}
\end{prop}

Since our studies have for origin the article of Chapuy \cite{chapuy} which proved (\ref{chap}) on random permutations, we give now the complete behavior of the above statistics in this case.
\begin{thm}
\label{permut}
Assume that $U$ is a random permutation of $[n]$. We have
\begin{eqnarray}
\label{chap}n^{-1/2}W\sn &\convlaw& W^{(\infty)}\\
\label{nchap1}
n^{-1/2}\left(\mathcal T\sn - \mathbb E  \mathcal T\sn\right) &\convlaw& B_{0,0}\\
\label{nchap2}n^{-1/2}\mathcal Z\sn &\convlaw& W^{(\infty)}\,.
\end{eqnarray}
\end{thm}

\begin{rem}
\label{remn}
Let $M^{p,q} = U^{p,q}\left(U^{p,q}\right)^*$ and $\mathcal M^{s,t} = \mathcal U^{s,t}\left(\mathcal U^{s,t}\right)^*$.
For $s,t$ fixed, the random variables $\mathcal T\sn_{s,t}$ and $T\sn_{s,t}$ are linear 
functionals of the empirical spectral distribution of $M^{\lfloor ns\rfloor, \lfloor nt\rfloor}$ and 
$\mathcal M^{s,t}$ respectively. For
 classical models in Random Matrix Theory, the convergence of fluctuations of such linear functionals do not need a normalizing factor, 
since the variance is bounded (the eigenvalues are repelling each other).  
Here, this is indeed the case for $T\sn_{s,t}$ (see \cite{DuPa} for the complete behavior for general tests functions). 
But, in the case of $\mathcal T\sn_{s,t}$, we have 
$ \Var \mathbb E [\mathcal T\sn_{s,t} | R_1, \dots, R_n, L_1, \dots, L_n] = O(n)$, which demands a normalization. 
\end{rem}

\begin{rem}
Going back to the decomposition (\ref{anova}), gathering the convergences in (\ref{convw}) (\ref{nchap1}) and (\ref{nchap2}) and the covariances,
we recover the identity
\[B_{0,0} \law W^{\infty} + \mathcal   W^{\infty}\]
where the two processes in the RHS are independent, fact which was quoted in \cite{DPY} section 2.
\end{rem}
\section{Proofs}
\subsection{Proof of Theorem \ref{onedim}}
It is known that the vector $(|U_{i1}|^2, i=1, \dots, n)$ is distributed uniformly on the simplex.
Actually, the result is a consequence of the following proposition. 
\begin{prop}
\label{onedimbeta}
Let $(\mathbf u_i, i=1,\dots,n)$ Dirichlet $(\beta', \dots, \beta')$ distributed. 
\begin{enumerate}
\item The process $n^{1/2}\{\sum_{i=1}^{\lfloor ns\rfloor} \left(\mathbf u_i - 1/n\right)\ , \  s\in [0,1]\}$ converges in distribution to $\sqrt{\beta'^{-1}} B_0$.
\item The process $n^{1/2}\{\sum_1^n \mathbf u_i \left(1_{R_i \leq s} - s\right)\ , \  s\in [0,1]\}$ converges in distribution to $\sqrt{1 + \beta'^{-1}}B_0$.
\end{enumerate}
\end{prop}
\proof
It is easy to see that
\begin{equation}
\label{idlawu}(\mathbf u_{1}, \cdots, \mathbf u_{n}) \law\left(\frac{ \mathbf g_1}{\mathbf g_1 + \cdots + \mathbf g_n}, \cdots , \frac{ \mathbf g_n}{\mathbf g_1 + \cdots + \mathbf g_n}\right)\end{equation}
with $\mathbf g_i$ independent and gamma distributed with parameter $\beta'$.

(1)  The above representation yields the following equality in law between processes
\[n^{1/2}\{\sum_{i=1}^{\lfloor ns\rfloor} \left(\mathbf u_i - 1/n\right)\law \frac{n}{\S_n}\left(\frac{\S_{\lfloor ns\rfloor}-\beta'\lfloor ns\rfloor}{\sqrt n}- \frac{\lfloor ns\rfloor}{n} \frac{\S_n - n \beta'}{\sqrt n}\right) \]
where $\S_k = \sum_{i=1}^k \mathbf g_i $. The WLLN and the Donsker's theorem gives 
the convergence to $  \sqrt{\beta'^{-1}} (B(s) -sB(1))$ where $B$ is the standard Brownian motion.

(2)  The process $n^{1/2}\{\sum_1^n \mathbf u_i \left(1_{R_i \leq s} - s\right)\ , \  s\in [0,1]\}$ is an example of a weighted empirical process. We may apply the Theorem 1.1  of Koul and Ossiander \cite{KoulOssiander} p.544. 
Set $\mathbf v_{ni} = n \mathbf u_i$ be the weights. Under the following conditions on the weights:
\begin{enumerate}
\item
$\left( n^{-1}\sum_i \mathbf v_{ni}^2\right)^{1/2}\convprob \gamma$ 
\item
$n^{-1/2}\max_i \mathbf v_{ni} \convprob 0$,
\end{enumerate}
 the process converges to the product of $\gamma$ by an independent Brownian bridge.
For the first one, we start from the representation (\ref{idlawu}), apply the SLLN twice and get
\[\left( n^{-1}\sum_i \mathbf v_{ni}^2\right)^{1/2} \convprob \frac{\left(\mathbb E g_1^2\right)^{1/2}}{\mathbb E g_1} = \sqrt{1 + \beta'^{-1}}\,.\]
For the second one, we have:
\[n^{-1/2}\max_i \mathbf v_{ni} = n^{1/2}\frac{\max \{\mathbf g_i , i=1, \dots, n\}}{\mathbf g_1 + \cdots + \mathbf g_n}.\]
Appealing again to the SLLN, it is enough to prove that $\max \{\mathbf g_i , i=1, \dots, n\} = o(n^{1/2})$ in probability, which is clear (this maximum is of order $\log n$). $\Box$

\subsection{Proof of Theorem \ref{permut}}
We notice that the double sums in the definitions of $\mathcal T\sn$ or $\mathcal Z\sn$ are actually single sums where $j = \pi(i)$ for $\pi$ a random permutation. Moreover, by exchangeability of the variables $C_j$, these single sums have (as processes in $s,t$) the same distributions as if $\pi$ was the identity permutation.  In this case, we may apply the results of \cite{Neu}, since 
$n^{-1/2}[\mathcal T\sn -\mathbb E \mathcal T\sn]$ and $n^{-1/2}\mathcal Z\sn$ are precisely $X_n^F$ and $Y_n$ therein, respectively.

\subsection{Proof of Proposition \ref{twodimprop}}
Under the DFT model, we have $\mathcal Z\sn = F\sn \otimes G\sn$, so we have only to study the Haar model, 
in $\mathbb U(n)$ and in $\mathbb O(n)$. As long as possible, we keep the notation $U_{ij}$ for the generic element of $U$ in both cases, but we use the notation $O_{ij}$ when we need to stress that the underlying unitary matrix is sampled according to the Haar measure in $\mathbb O (n)$.  

First, we prove (\ref{maincv}), i.e. the convergence of finite dimensional distributions, beginning with the 1-dimensional marginals since the notations are easier to follow for the reader, and then using the Cramer-Wold device. In a last part, we prove the tightness, using a criterion due  to Davydov and Zitikis \cite{DaZ}.

\subsubsection{1-Marginals}

We first consider the convergence of the 1-marginals, i.e. for fixed  $s,t \in [0,1]$, we show the convergence of $\mathcal Z^{(n)}(s,t)$ to $\sqrt{\beta'^{-1}} W^{(\infty)} (s,t) + B_0 \otimes B_0 (s,t)$. 
By a scaling property, this is equivalent to prove the weak convergence of $\displaystyle \frac{\mathcal Z^{(n)}(s,t)}{\sqrt{s(1-s)t(1-t)}}$ towards $\sqrt{\beta'^{-1}}N_1 + N_2N_3$ where $N_i$ are independent standard Gaussian variables.  According to the next Proposition, whose proof is given in the Appendix,  we can replace  the independent centered random variables 
\[R^i := \frac{{\mathbf 1}_{R_i \leq s}-s}{\sqrt{s(1-s)}} \ , \  C^j :=  \frac{{\mathbf 1}_{C_j \leq t}-t}{\sqrt{t(1-t)}} \  ,  \  i,j =1, \dots, n \]
by  independent $\mathcal N(0,1)$ random variables $X_i$ and $Y_j$. 
\begin{prop}
\label{Lin1}
Let $(R_i)$ and $(C_j)$ two independent sequences of iid centered variables, with variance 1, and finite third moment. Consider also $(X_i)$ and $(Y_j)$ two independent sequences of iid standard Gaussian  variables.
We define $A_n = \sum_{i,j=1}^n |U_{ij}|^2 R_i C_j$ and $B_n = \sum_{i,j=1}^n |U_{ij}|^2 X_i Y_j$.
Then, $A_n$ and $B_n$ have the same limiting distribution.
\end{prop}
\medskip

We thus  study the bilinear non symmetric form
\[\Sigma_n := \sum_{i,j=1}^n X_i |U_{ij}|^2 Y_j\]
built from the non symmetric matrix $\tilde{U}= \left(|U_{ij}|^2\right)_{i,j \leq n}$. This matrix  
is a so-called "doubly stochastic" matrix (\cite{unistoch}, \cite{chaf}).
As a Markovian  matrix, it has $1$ as dominant eigenvalue.
\begin{prop}
\label{propSigma}
The sequence of random variables
\[\Sigma_n := \sum_{i,j=1}^n X_i |U_{ij}|^2 Y_j\]
converges in distribution towards $\sqrt{\beta'^{-1}} N_1 + N_2 N_3$ where $N_1, N_2,N_3$ are independent and $\mathcal N (0,1)$ distributed.
\end{prop}

\proof
We know that $\mathbb E |U_{ij}|^2 = n^{-1}$. We may remark that 
\[\mathbb E [\Sigma_n |X,Y] =  n^{-1}\sum_{i,j=1}^n X_i Y_j = \left(\frac{X_1+ \dots+X_n}{\sqrt n}\right) \left(\frac{Y_1+ \dots+Y_n}{\sqrt n}\right)=: S'_n\]
and that  $S'_n \law N_2N_3$. Set 
\[S_n := \Sigma_n  - S'_n =  \sum_{i,j=1}^n X_i V_{ij}Y_j\] 
with \[ V_{ij}= |U_{ij}|^2 - n^{-1}\,.\]
We will use the characteristic function and conditioning. 
Set $(X,V) := (X_i, i \leq n , V_{ij}, i, j \leq n)$.  Conditionnally upon $(X, V)$,  the vector $(S_n , S'_n)$ is Gaussian and 
\be\cov(S_n, S'_n | X, V) &=& n^{-1}\sum_{i,j,k,\ell=1}^n X_iX_kV_{k\ell} \cov(Y_j, Y_\ell)\\
&=&  n^{-1}\sum_{i,j,k=1}^n X_iX_kV_{kj} =  n^{-1}\sum_{i,k=1}^n X_iX_k \left(\sum_{j=1}^n V_{kj}\right) = 0\ee
since $U$ is unitary. The variables $S_n$ and $S'_n$ are then  independent, conditionnally upon $(X,V)$. 
Moreover
\ben
\mathbb E [\exp\ii \theta S'_n | (X,V)] = \mathbb E [\exp\ii \theta S'_n | X ] = \exp - \frac{\theta^2}{2}\frac{(\sum_{i=1}^n X_i)^2}{n}
\een
and
\ben
\mathbb E [\exp \ii \theta S_n |(X,V)] = \exp -\frac{\theta^2}{2} \widehat S_n
\een
where 
\[\widehat S_n := \sum_{j=1}^n \left(\sum_{i=1}^n V_{ij} X_i\right)^2 =  \sum_{i,j=1}^n X_iX_j H_{ij}\,,\]
with
\[H_{ij} = (VV^*)_{ij} \ , \ i,j =1, \dots, n\,.\]
This implies
\ben
\label{crucial}
\mathbb E \exp{\ii \theta\Sigma_n} =\mathbb E \exp -\frac{\theta^2}{2}\left(n^{-1}(\sum_{i=1}^n X_i)^2 +  \widehat S_n \right)
\een
Since $(\sum X_i)^2 /n$ is distributed as $N^2$ (where $N$ is $\mathcal N(0,1)$), if we prove that 
\ben
\label{cj2}
\widehat S_n :=  \sum_{j=1}^n \left(\sum_{i=1}^n V_{ij} X_i\right)^2 \convprob \beta'^{-1},
\een
 we will conclude that
\[\mathbb E \exp \ii \theta \Sigma_n \rightarrow \mathbb E \exp -\frac{\theta^2}{2} (N^2 +\beta'^{-1})\]
and this limit is the characteristic function of $N_2N_3 + \sqrt{\beta'^{-1}} N_1$.

  It is clear that  (\ref{cj2}) is implied by
\ben
\label{cj3e}
\lim_n \mathbb E \widehat S_n &=& \beta'^{-1}\\
\label{cj3v}
\lim_n  \mathbb E (\widehat S_n)^2 &=& \beta'^{-2}\,.
\een
To prove these assertions, we will need some joint moments of elements of the matrix $H$, which are themselves affine functions of moments of the matrix $U$.

The first limit in (\ref{cj3e}) is easy to obtain since
\ben
\nonumber
\mathbb E \widehat S_n &=& \sum_i  \mathbb E H_{ii} = n \mathbb E H_{11} = n^2 \mathbb E (V_{11})^2 = n^2 \Var V_{11}\\
\label{cj3er}
&=&  \begin{cases}n^2  \Var |O_{11}|^2 = \frac{2(n+3)}{n+2}\;\;&\mbox{in the orthogonal case} ,\\
n^2  \Var |U_{11}|^2 = \frac{n-1}{n+1}\;\;&\mbox{in the unitary case,}
\end{cases} \een
where we refer to \eqref{momentX}.

To prove  (\ref{cj3v}) we expand $\widehat S_n^2$:
\[ \widehat S_n^2 = \sum_{i,j,k,\ell=1}^n X_iX_jX_kX_\ell H_{ij}H_{k\ell}\,.\]
 In the expectation, the only non vanishing terms are obtained according to the decompostion 
\begin{eqnarray}
\label{devhat}
\mathbb E \widehat S_n^2 &=& \sum_{i=j, k=\ell\not= i} +  \sum_{i=k, j=\ell\not= i} +  \sum_{i=\ell, j=k\not= i}+ \sum_{i=j=k=\ell}\\
&=& (1)+ (2) + (3) + (4)\,.\end{eqnarray}
Since $ (\mathbb E X_1^2)^2 = 1$ and $\mathbb E X_1^4 = 3$ we have successively 
\ben
\label{dewhat1}
(1) &=& 
  n(n-1)  \mathbb  E\left(H_{11}H_{22}\right)\\
\label{dewhat2}
(2) = (3) &=&    n(n-1)  \mathbb E (H_{12}^2)\\
\label{devhat4}
(4) &=&  3n  \mathbb E H_{11}^2\,.
\een
Now
\ben
\label{h1}
\mathbb E H_{11}^2 &=& n \mathbb E V_{11}^4 + n(n-1) \mathbb E \left(V_{11}^2V_{12}^2\right)\\
\label{h2}
 \mathbb  E\left(H_{11}H_{22}\right) &=& n \mathbb E \left(V_{11}^2 V_{21}^2\right) + n(n-1) \mathbb E \left(V_{11}^2 V_{22}^2\right)\\
\label{h3}
\mathbb E (H_{12}^2)  &=& n\mathbb E\left(V_{11}^2 V_{21}^2 \right) + n(n-1) \mathbb E \left(V_{11}V_{12}V_{21}V_{22}\right)\,.
\een
this last equality coming from the symmetry of $H$.

From Lemma  \ref{lemme1} (see \eqref{Weing-unit1} with $k=4$), 
\begin{equation} \label{S_n}
 \E \widehat S_n^2 = n^2(n-1)^2 \mathbb E \left(V_{11}^2 V_{22}^2\right) + 2 n^2(n-1)^2
\mathbb E \left(V_{11}V_{12}V_{21}V_{22}\right) + o(1). 
\end{equation}
 and from Lemma \ref{lem1}, we conclude that (\ref{cj3v}) holds, which ends the proof of \ref{propSigma}. $\Box$

\subsection{Finite-dimensional marginals}
We now consider the convergence of the finite dimensional distributions, following the same scheme of proof as in the case of the 1-dimensional marginal. The Lindeberg's strategy statement is now:
\begin{prop} 
\label{Lin2}
Let $(\beta_i(s))$ and $(\gamma_i(s))$ be two independent sequences of independent Brownian bridges. 
Define $G_n(s,t) = \sum_{i,j=1}^n |U_{ij}|^2 \beta_i(s) \gamma_j(t)$. Then, the processes $\mathcal Z_n$ and $G_n$ have the same finite dimensional  limiting distributions.
\end{prop}

We then consider the convergence of the finite-dimensional distributions of the process
$$G_n(s,t) =  \sum_{i,j} |U_{ij}|^2 \beta_i(s) \gamma_j(t) = S_n(s,t) + S'_n(s,t)$$ 
where 
$ S'_n(s,t) = \frac{1}{n} \sum_{i,j}  \beta_i(s) \gamma_j(t) $,
$S_n(s,t) = \sum_{i,j}V_{ij} \beta_i(s) \gamma_j(t)$ and $(\beta_i)$, $(\gamma_j)$ are two independent sequences of independent bridges. \\
Let $(s_1,t_1), \ldots, (s_K,t_K) \in [0,1]^2$, $\alpha_1, \ldots, \alpha_K \in \mathbb R$ and define
$$\Sigma_n = \sum_{l=1}^K \alpha_l G_n(s_l,t_l) : = S_n + S'_n$$
according to the decomposition of $G_n$.
Conditionnally to $(V, \beta_i, i \leq n)$, $(S_n,S'_n)$ is a Gaussian vector with covariance
\begin{eqnarray*}
\cov(S,S'| (V, \beta)) &=&  \frac{1}{n} \sum_{l,l'=1}^K \sum_{i,j,i',j'= 1}^n \alpha_l \alpha_{l'}V_{ij} \beta_i(s_l) \beta_{i'}(s_{l'} )\cov(\gamma_j(t_l), \gamma_{j'}( t_{l'})) \\
&=& \frac{1}{n}  \sum_{l,l'=1}^K  \sum_{i,i'=1}^n\alpha_l \alpha_{l'}  \beta_i(s_l) \beta_{i'}(s_{l'}) g(t_l,t_{l'}) \sum_{j=1}^n V_{ij} \\
&=& 0
\end{eqnarray*}
where $g(t,t')$ denotes the covariance of the Brownian bridge and we use that $\sum_{j} V_{ij} = 0$.
Thus $S_n$ and $S'_n$ are  conditionnally independent given $(V,\beta)$. Thus,
\be
  \mathbb E\left[ \exp{\ii \theta\Sigma_n}  | (V, \beta)\right]
&=&  \mathbb E \left( \mathbb E\left[ \exp{\ii \theta S'_n}  | (V, \beta)\right] \mathbb E\left[ \exp{\ii \theta S_n}  | (V, \beta)\right]\right)\\
\ee
\be 
\mathbb E\left[ \exp{\ii \theta S'_n}  | (V, \beta)\right] &= &\exp \left(- \frac{\theta^2}{2n} \left( \sum_{l,l'=1}^K \alpha_l \alpha_{l'} g(t_l,t_{l'}) (\sum_{i=1}^n \beta_i(s_l)) (\sum_{i=1}^n \beta_i(s_{l'})\right) \right) \\
&\law&  \exp \left(- \frac{\theta^2}{2} \left( \sum_{l,l'=1}^K \alpha_l \alpha_{l'} g(t_l,t_{l'}) \beta(s_l) \beta(s_{l'})\right) \right)\,,
\ee
and
\be
E\left[ \exp{\ii \theta S_n}  | (V, \beta)\right]  =  \exp \left(- \frac{\theta^2}{2} \widehat{S}_n\right)
\ee
where 
 $$\widehat{S}_n = \sum_{l,l'=1}^K \alpha_l \alpha_{l'} g(t_l,t_{l'}) \sum_{i,k=1}^n H_{i,k} \beta_i(s_l) \beta_k(s_{l'}).$$
 We now prove the convergence in probability of $\widehat{S}_n$ to 
\[L := \beta'^{-1} \sum_{l,l'=1}^K \alpha_l \alpha_{l'} g(s_l, s_{l'}) g(t_l, t_{l'})\,.\]
We have first $$\mathbb E(\widehat{S}_n) = 
\sum_{l,l'=1}^K \alpha_l \alpha_{l'} g(t_l,t_{l'})g(s_l,s_{l'}) 
\sum_{i=1}^n \E(H_{ii})$$
and  $\sum_{i} \E(H_{ii}) = n \E(H_{11}) \longrightarrow \beta'^{-1}$. Now, 
\begin{eqnarray*}
\widehat{S}_n^2 &=& \sum_{l,l',p,p'=1}^K  \sum_{i,j,k,q=1}^n \alpha_l \alpha_{l'} \alpha_p \alpha_{p'} g(t_l,t_{l'}) g(t_p, t_{p'}) \\
& &  \quad H_{ik} H_{jq} \beta_i(s_l) \beta_{k} (s_{l'}) \beta_{j}(s_p) \beta_q(s_{p'})\,.
\end{eqnarray*}
 Taking the expectation, non zero terms are obtained when the indexes $i,j,k,q$ are equal 2 by 2.
 Moreover, the only non null contribution at the limit is given for the combination $i=k, j=q$, $i\not= j$ from Lemma  \ref{lem1}. We thus obtain that $\E(\widehat{S}_n^2 )$ converges to 
  $$ \beta'^{-2} \sum_{l,l',p,p'=1}^K  \alpha_l \alpha_{l'} \alpha_p \alpha_{p'} g(t_l,t_{l'}) g(t_p, t_{p'}) 
  g(s_l,s_{l'}) g(s_p,s_{p'}) = L^2$$
proving the convergence in probability of  $\widehat{S}_n$ as desired.
 Thus, we have obtained 
 $$ \E(\exp{\ii \theta \Sigma_n})  \longrightarrow  \E\left( \exp - \frac{\theta^2}{2} \left(\sum_{l,l'} \alpha_l \alpha_{l'} g(t_l,t_{l'}) \beta(s_l) \beta(s_{l'}) + L
\right) \right)\,.$$
Since 
\[L = \mathbb E \left(\sum_l \alpha_l \sqrt{ \beta'^{-1}}W^{(\infty)}(s_l, t_l)\right)^2\]
and
\[\mathbb E\left( \exp -\frac{\theta^2}{2} \left(\sum_{l,l'} \alpha_l \alpha_{l'} g(t_l,t_{l'}) \beta(s_l) \beta(s_{l'})\right)\right) = \mathbb E \left(\exp i\theta  \sum_l \alpha_l B^0\otimes B^0 (s_l, t_l)\right),\]
we have proved \eqref{maincv} in the sense of the finite dimensional distributions.

\subsection{Tightness of $n^{-1/2}\mathcal Z\sn$}
For the tightness, owing to the structure of the process $\mathcal Z_n$ we apply a criterion 
due to Davydov and Zitikis \cite{DaZ} which can be reformulated as follows.
\begin{thm}[Davydov-Zitikis]
\label{DaZi}
Let $\xi\sn$ be a sequence of stochastic processes in $D([0,1)^2$ such that $\xi\sn \convfidi \xi$ and $\xi\sn(0) = 0$. Assume that
\begin{enumerate}
\item there are constants $\alpha\geq \beta > 2, c \in (0, \infty)$ and $a_n \downarrow 0$ such that 
for all $n \geq 1$, we have 
\ben
\label{crit1}
\sup_{s', t'} \mathbb E \left(|\xi\sn (s+s', t+t') - \xi\sn (s', t')|^\alpha|\right)\leq  c ||(s,t)||^\beta\,,
\een
whenever $||(s,t)|| \geq a_n$, where $||(s,t)|| = \max(|s|, |t|)$.
\item the process $\xi\sn$ can be written as the difference of two coordinate-wise non-decreasing processes $\xi\sn_1$ and $\xi\sn_2$ such that $\xi\sn_2$ satisfies 
\ben
\nonumber
\sup_{s,t \leq 1}\max \{ \xi\sn_2\left(s, t + a_n\right)- \xi\sn_2\left(s,t\right) , \xi\sn_2\left(s+a_n, t\right)- \xi\sn_2\left(s,t\right)\}=o_P(1)\,,\\
\label{crit2}
\een 
\end{enumerate}
then the process $\xi\sn$ converges weakly to 
$\xi$. 
\end{thm}
We will prove the two following lemmas.

\begin{lem}\label{lcrit1} $\mathcal Z\sn$ satisfies (\ref{crit1}).
\end{lem}
\begin{lem}For any sequence $(b_n)$ with $b_n \rightarrow + \infty$, 
\label{lcrit2}
 $b_n^{-1} \mathcal Z\sn$ satisfies (\ref{crit2}).
\end{lem}

\noindent{\sl{Proof  of Lemma \ref{lcrit1}}}.
Let us compute the  moment of order 6 of the increments. Let $s, s', t,t' \in [0,1]$ such that $s+s' \leq 1, t+t' \leq 1$ and $B = ]s', s+s'] \times ]t', t+t']$. Recall that the increment of $\mathcal Z_n$ over $B$ is given by 
$$\mathcal Z_n(B) = \mathcal Z_n(s+s', t+t') - \mathcal Z_n(s+s', t') - \mathcal Z_n(s', t+t') + \mathcal Z_n(s', t').$$
It is easy to see that 
$$ \E( \mathcal Z_n(B)^6) =  \E( \mathcal Z_n(s,t)^6).$$
We have, using the notation $B_i= (1_{(L_i \leq s)} -s)$ and $B_j'= (1_{(C_j \leq t)} -t)$ :
\begin{equation}
\E (\mathcal Z_n(B)^6) = \sum_{i_k, j_k, k=1, \ldots 6} \E(\prod_{k=1}^6 |U_{i_kj_k}|^2) \E(\prod_{k=1}^6 B_{i_k}) \E(\prod_{k=1}^6 B'_{j_k}).
\end{equation}
Since the $B_i$ and $B'_j$ are independent and centered, in the RHS of the above equation, the non null term in the sum are obtained when the $i_k$ (resp. the $j_k$) are equal at least 2 by 2.
We now use the following properties: for some $C>0$,
\begin{itemize}
\item $ |\E(B_i^k)| \leq Cs$ for $2\leq k\leq 6$
\item $ |\E((B'_j)^k)| \leq Ct$ for $2\leq k\leq 6$
\item $\E(\prod_{k=1}^6 |U_{i_kj_k}|^2) = O(\frac{1}{n^6})$ (see \eqref{Weing-unit1})
\end{itemize}
It follows that ($C$ may change from line to line)
\ben \label{tight}
\E (\mathcal Z_n(B)^6) \leq C \left(\frac{st}{n^4} +  \frac{st^2}{n^3} + \frac{st^3}{n^2} + \frac{s^2t}{n^3}+ \frac{s^2t^2}{n^2}+ \frac{s^2t^3}{n}+ s^3t^3\right).
\een
Note that if $s \geq \frac{1}{n}$ and $t \geq \frac{1}{n}$, 
\ben \label{tight1}
\E (\mathcal Z_n(B)^6) \leq C  s^3t^3.
\een
Define $||(s,t)|| = \sup(|s|, |t|)$. If $||(s,t)|| \geq \frac{1}{n}$, we can see from \eqref{tight} that 
\ben \label{tight2}
\E(\mathcal Z_n(B)^6) \leq C  ||(s,t)||^3.
\een
We now consider the increment $\mathcal Z_n(s+s', t+t') - \mathcal Z_n(s', t')$. It is easy to see that 
$$ \mathcal Z_n(s+s', t+t') - \mathcal Z_n(s', t') = \mathcal Z_n(B) + \mathcal Z_n(B_1) + \mathcal Z_n(B_2)$$
where $B_1 = [0,s'] \times [t', t'+t]$ and $B_2 = [s', s'+s] \times [0,t']$.
Therefore,
$$\E(|\mathcal Z_n(s+s', t+t') - \mathcal Z_n(s', t')|^6) \leq C (\E (\mathcal Z_n(B)^6) + \E (\mathcal Z_n(B_1)^6) + \E (\mathcal Z_n(B_2)^6)).$$
From \eqref{tight}, we can see that in any case, for $i=1,2$, if $||(s,t)|| \geq \frac{1}{n}$ :
$$ \E (\mathcal Z_n(B_i)^6) \leq C   ||(s,t)||^3.$$
For example, if $ s \leq \frac{1}{n} \leq t$,
$$  \E (\mathcal Z_n(B_2)^6)) \leq C \frac{s}{n^2} \leq C st^2 \leq C t^3.$$
We have thus proved the following: if $||(s,t)|| \geq \frac{1}{n}$, 
$$\E(|\mathcal Z_n(s+s', t+t') - \mathcal Z_n(s', t')|^6) \leq C ||(s,t)||^3\,,$$
which is exactly (\ref{crit1}) with $a_n = n^{-1}$, $\alpha = 6$ and $\beta= 3$.  $\Box$
\bigskip

\noindent{\sl{Proof  of Lemma \ref{lcrit2}}}.
Indeed, we can write:
\[\mathcal Z_n (s,t) = \Xi_1 (s,t) - \Xi_2(s,t)\]
with
\[\Xi_1(s,t) = \sum_{i,j}|U_{ij}|^2 1_{L_i \leq s}1_{C_j \leq t} +st  \  , \  
\Xi_2 (s,t) = s \sum_j 1_{C_j \leq t} + t \sum_i 1_{L_i \leq s}\]
and both these processes are coordinate-wise non-decreasing.
Let us now check (\ref{crit2}) with $\xi\sn_1 = \Xi_1/b_n$ and $\xi\sn_2 = \Xi_2/b_n$.
We have to prove that for $\delta = 1/n$,
\ben
\nonumber
\sup_{s,t \leq 1}\max \{ \Xi_2\left(s, t + \delta\right)- \Xi_2\left(s,t\right) , \Xi_2\left(s+\delta, t\right)- \Xi_2\left(s,t\right)\}=o_P(b_n)\,.\\
\label{cor12}
\een
Owing to the symmetry of the roles played by $C_j$ and $L_i$, we will focus on the first term in the above maximum.
We have
\ben\nonumber\Xi_2\left(s, t + \delta\right)- \Xi_2\left(s,t\right) &=& s \sum_j {\mathbf 1}_{t\leq C_j\leq t+\delta} + \delta \sum_i {\mathbf 1}_{L_i \leq s}\\
&\leq& \sum_j {\mathbf 1}_{t\leq C_j\leq t+\delta} +  n\delta\\
&:=& N_n(t, 1/n) + 1.\een
We have to prove that, for every $\varepsilon > 0$,
\ben\label{Davy}\mathbb P (\sup_t  N_n(t, 1/n) > \varepsilon b_n) \rightarrow 0\,.\een
or, equivalently,
\ben\label{Davyd}\mathbb P (\sup_t  N_n(t, 1/n) > \lfloor\varepsilon b_n\rfloor) \rightarrow 0\,.\een
The event $\{\sup_t  N_n(t, 1/n) > \lfloor\varepsilon b_n\rfloor\}$ means that there is a subinterval of $[0,1]$ of length $1/n$ which contains at least $\lfloor\varepsilon b_n\rfloor$ points of the sample $(C_1, \dots, C_n)$.
Denote by $\left(C_{(1)}\sn \leq C_{(2)}\sn, \dots, \leq C_{(n)}\sn\right)
$ the reordered sample. We have (denoting $C_0 = 0$ by convention)
\ben
\{\sup_t  N_n(t, 1/n) > \lfloor\varepsilon b_n\rfloor\} \subset \{\exists\ 0 \leq  k \leq n : C_{\left(k+ \lfloor b_n \varepsilon\rfloor\right)}\sn - C_{(k)}\sn \leq 1/n\}
\een
hence, by the union bound
\ben\label{ub}\mathbb P  (\sup_t  N_n(t, 1/n) > \lfloor\varepsilon b_n\rfloor) \leq \sum_k \mathbb P\left( C_{\left(k+ \lfloor b_n\varepsilon\rfloor\right)}\sn - C_{(k)}\sn \leq 1/n\right)\,.
\een
It is well known that the spacings  follow the Dirichlet distribution of parameter $(1, 1, \dots, 1)$ so that
\ben C_{(k+r)} - C_{(k)} \law \frac{\mathbf g_1 + \dots + \mathbf g_r}{\mathbf g_1 + \dots + \mathbf g_n}\,,\een
where $\mathbf g_i$ are i.i.d. and exponential.
From (\ref{ub}), we get
\ben
\mathbb P  (\sup_t  N_n(t, 1/n) > \lfloor\varepsilon b_n \rfloor) \leq n \mathbb P \left(\frac{\mathbf g_1 + \dots + \mathbf g_{\lfloor b_n\varepsilon\rfloor}}{\mathbf g_1 + \dots + \mathbf g_n} \leq \frac{1}{n}\right)\leq \\
\nonumber 
\leq  n \mathbb P(\mathbf g_1 + \dots + \mathbf g_n > 2n) + n  \mathbb P(\mathbf g_1 + \dots + \mathbf g_{\lfloor b_n\varepsilon\rfloor} \leq 2)\,.
\een
But we know (Chernov bound) that 
\ben
\mathbb P (\mathbf g_1 + \dots + \mathbf g_n > 2n) &\leq& \exp - n h(2)\\
\mathbb P (\mathbf g_1 + \dots + \mathbf g_{\lfloor b_n\varepsilon\rfloor} \leq 2)&\leq&  \exp - \lfloor b_n\varepsilon\rfloor h\left(\frac{2}{ \lfloor b_n\varepsilon\rfloor}\right)
\een
where
\[h(x) = x- 1 - \log x\,.\]
We conclude that (\ref{crit2}) is fulfilled.  $\Box$

\section{Moments}
\label{prelimunit}

If $U$ is Haar distributed matrix on $\mathbb U(n)$ or $O$ Haar distributed on $\mathbb O(n)$,  let us denote by $\mathbf u$ the generic element $|U_{ij}|^2$ or $O_{ij}^2$. We know that $\mathbf u$ follows the beta distribution on $[0,1]$ with parameter $(\beta', (n-1)\beta')$ so that 

\begin{eqnarray}
\label{momentX}
\mathbb E |U_{i,j}|^2 = \frac{1}{n} \ , \ \mathbb E |U_{i,j}|^4 = \frac{2}{n(n+1)} \ , \ \Var |U_{i,j}|^2 =  \frac{n-1}{n^2(n+1)}\\
\mathbb E |O_{i,j}|^2 = \frac{1}{n} \ , \ \mathbb E |O_{i,j}|^4 = \frac{3}{n(n+2)} \ , \ \Var |O_{i,j}|^2 =  \frac{2(n+3)}{n^2(n+2)}\,, 
\end{eqnarray}
and more generally
\begin{eqnarray} \label{A1}
\mathbb E |U_{ij}|^{2k} &=& \frac{(n-1)! k!}{(n-1+k)!}\\
 \label{A1o}
\mathbb E |O_{ij}|^{2k} &=& \frac{(2k)!!}{n(n+2) \dots (n+k-2)}\,.
\end{eqnarray}

\begin{lem} \label{lemme1} Let $k \in \NN$,
\begin{enumerate}
\item  For every choice of indices ${\mathbf i} = (i_1, \ldots, i_{k})$ and 
 ${\mathbf j} = (j_1, \ldots, j_{k})$,
\begin{equation} \label{Weing-unit1}
\mathbb E \left(|U_{i_1 j_1}|^2 \ldots |U_{i_k j_k}|^2 \right) = O(\frac{1}{n^k})
\end{equation}
\item  For every choice of indices ${\mathbf i} = (i_1, \ldots, i_{k})$ and 
 ${\mathbf j} = (j_1, \ldots, j_{k})$,
\begin{equation} \label{Weing-unit1}
\mathbb E \left(|O_{i_1 j_1}|^2 \ldots |O_{i_k j_k}|^2 \right) = O(\frac{1}{n^k})
\end{equation}
\item For every choice of indices ${\mathbf i} = (i_1, \ldots, i_{k})$ and 
 ${\mathbf j} = (j_1, \ldots, j_{k})$, for the unitary and orthogonal cases,
\begin{equation} \label{Weing-unit1}
\mathbb E \left(V_{i_1 j_1} \ldots V_{i_k j_k} \right) = O(\frac{1}{n^k})
\end{equation}
\end{enumerate}
\end{lem}
\proof

(1) follows from the Weingarten formula giving the moments of Haar unitary matrix coefficients (see \cite{CoSn}, Corollary 2.4).
\begin{equation} \label{Weing-unit2}
\mathbb E\left(  U_{i_1j_1} \ldots U_{i_k, j_k} \bar U_{i_{\bar 1}j_{\bar 1}}\ldots \bar U_{i_{\bar k}, j_{\bar k}}\right) 
= \sum_{\alpha, \beta \in {\mathcal S}_{k}} \tilde\delta_{\mathbf i}^{\alpha}  \tilde\delta_{\mathbf j}^{\beta} \W(n , \beta\alpha^{-1})
\end{equation}
where $\tilde\delta_{\mathbf i}^{\alpha} = 1$ if $i(s)= i(\overline{\alpha(s)})$ for every $s \leq k$ and $0$ otherwise, 
and the asymptotics of the Weingarten function (see \cite{CoSn}, Proposition 2.6):
\begin{equation} \label{asympW}
\W(n ,\sigma) = O(n^{-k - |\sigma|})
\end{equation}
 where $|\sigma| = k - \#(\sigma)$. Note that the maximal order is obtained for the identity permutation $Id$ for which $|Id| = 0$ and in this case, $\W(n ,Id) = n^{-k} (1+ o(1))$.

(2) follows first from the Weingarten formula giving the moments of Haar orthogonal matrix coefficients (see \cite{CoSn}, Corollary 3.4).
For every choice of indices ${\mathbf i} = (i_1, \ldots, i_{k}, i_{\bar 1}, \ldots, i_{\bar k})$ and  
${\mathbf j} = (j_1, \ldots, j_{k}, j_{\bar 1}, \ldots, j_{\bar k})$,
\begin{equation} \label{Weing-orth}
\mathbb E \left(  O_{i_1j_1} \ldots O_{i_{k}j_{k}}O_{i_{\bar 1}j_{\bar 1}} \ldots O_{i_{\bar k}j_{\bar k}}\right)
 = \sum_{p_1, p_2 \in {\mathcal M}_{2k}} \delta_{\mathbf i}^{p_1}  \delta_{\mathbf j}^{p_2} \W^{{\mathbb O}(n)}(p_1, p_2)
\end{equation}
where ${\mathcal M}_{2k}$ denotes the set of pairings of $[2k]$, $\W^{{\mathbb O}(n)}$ is the orthogonal  Weingarten matrix  and
 $\delta_{\mathbf i}^{p_1}$ (resp. $\delta_{\mathbf j}^{p_2})$ is equal to $1$  or $0$ if $\mathbf i$ (resp. $\mathbf j$) is constant on each pair of $p_1$ (resp. $p_2$) or not. \\
 We then use asymptotics for the orthogonal Weingarten matrix (see \cite{CoSn}, Theorem 3.13):
 $$\W^{{\mathbb O}(n)}(p_1, p_2) = O(n^{-k-l(p_1,p_2)})$$
 for some metric $l$ on ${\mathcal M}_{2k}$.

(3)  follows from (1), resp. (2),  and the definition of $V$ in terms of $|U|^2$. $\Box$. \\

We now need to have precise asymptotics for some   moments of $U$ up to order 8.

\begin{prop}
\label{propu}
\begin{enumerate}
\item
\begin{equation} \label{1}
\E( |U_{11}|^4 |U_{22}|^4) = \frac{4}{n^4} + O(\frac{1}{n^5})
\end{equation}
\begin{equation} \label{2}
\E( |U_{11}|^2 |U_{12}|^2|U_{21}|^2 |U_{22}|^2) = \frac{1}{n^4} + O(\frac{1}{n^5})
\end{equation}
\begin{equation} \label{3}
\E( |U_{11}|^2 |U_{22}|^4) = \frac{2}{n^3} + O(\frac{1}{n^4})
\end{equation}
\begin{equation} \label{4}
\E( |U_{11}|^2 |U_{12}|^2 |U_{22}|^2) = \frac{1}{n^3} + O(\frac{1}{n^4})
\end{equation}
\item
\begin{equation} \label{1o}
\E( |O_{11}|^4 |O_{22}|^4) = \frac{9}{n^4} + O(\frac{1}{n^5})
\end{equation}
\begin{equation} \label{2o}
\E( |O_{11}|^2 |O_{12}|^2|O_{21}|^2 |O_{22}|^2) = \frac{1}{n^4} + O(\frac{1}{n^5})
\end{equation}
\begin{equation} \label{3o}
\E( |O_{11}|^2 |O_{22}|^4) = \frac{3}{n^3} + O(\frac{1}{n^4})
\end{equation}
\begin{equation} \label{4o}
\E( |O_{11}|^2 |O_{12}|^2 |O_{22}|^2) = \frac{1}{n^3} + O(\frac{1}{n^4})
\end{equation}
\end{enumerate}
\end{prop}
\proof  
Let us first prove (1).
We shall give the details for the first formula \eqref{1}, the other ones are similar.
From (\ref{Weing-unit2}) and \eqref{asympW}, the main contribution in (\ref{Weing-unit2}) is given by the pairs  $(\alpha, \beta)$ for which $\beta\alpha^{-1} = Id$.
We now give all the permutations $\alpha$ and $\beta$ giving a non null contribution in (\ref{Weing-unit2}) for the computation of $\E( |U_{11}|^4 |U_{22}|^4)$.
Tha admissible permutations $\alpha$ are given by their cycle decomposition:
$$\begin{array}{l}
\alpha_1 = Id = (1), (2), (3), (4) \\
 \alpha_2 = (12), (3),(4) \\
 \alpha_3 = (1),(2), (34)\\
 \alpha_4 = (12), (34)\\ 
 \end{array} $$
 and the same for the $\beta$. \\
 There are 4 pairs of permutations giving $\beta\alpha^{-1} = Id$, those corresponding to $\alpha = \beta$. Therefore, we obtain \eqref{1}. \\
 For \eqref{2}, the admissible $\alpha$ are the same and the corresponding $\beta$ are given by
 $$\begin{array}{l}
\beta_1 = Id = (1), (2), (3), (4) \\
\beta_2 = (13), (2),(4) \\
\beta_3 = (1),(3), (24)\\
\beta_4 = (13), (24).
\\ \end{array} $$
 Thus, there is only one pair $(\alpha_1, \beta_1)$ giving the main contribution.\\
 The proof for the moments of order 6 is similar. 

Let us now prove (2). We refer to  the paper \cite{BaCoS}. In particular, these authors define
\ben
\label{defI}
 I_n\begin{pmatrix} a&c\\b&d\end{pmatrix} = \mathbb E (O_{11}^a O_{12}^b O_{21}^c O_{22}^d)\ , \ 
I_{n-1}(a, b) = \mathbb E (O_{11}^aO_{12}^b).
\een
In the sequel we denote
\[m!! = (m-1) (m-3) \dots\]
We need
\[I_n\begin{pmatrix} 4&0\\0&4\end{pmatrix}, I_n\begin{pmatrix} 2&2\\2&2\end{pmatrix}, 	I_n\begin{pmatrix} 2&0\\0&4\end{pmatrix}, I_n\begin{pmatrix} 2&2\\0&2\end{pmatrix}\]
The first and the third are ruled by Theorem C therein, so that
\[I_n\begin{pmatrix} 2&0\\0&4\end{pmatrix} = \frac{n!  2!! 4!! (n+4)!!}{n!! (n+2)!! (n+5)!!} = \frac{3(n+3)}{n(n-1)(n+2)(n+4)}= \frac{3}{n^3}+ O(\frac{1}{n^4})\]
and
\begin{eqnarray*}
I_n\begin{pmatrix} 4&0\\0&4\end{pmatrix}&=& \frac{n!  4!! 4!! (n+6)!!}{(n+2)!! (n+2)!! (n+7)!!} \\&=&\frac{9 (n+3)(n+5)}{(n+6)(n+4)(n+2)(n+1) n(n-1) } 
= \frac{9}{n^4} + O(\frac{1}{n^5})
\end{eqnarray*}
Now, from Definition 3.1 therein
\[ I_n\begin{pmatrix} 2&2\\2&2\end{pmatrix} = [I_{n-1} (2,2)]^2\Phi\begin{pmatrix}2&2\\2&2\end{pmatrix}\]
and from Theorem 6.2
\[\Phi\begin{pmatrix}2&2\\2&2\end{pmatrix} = \frac{n-1)!!}{(n-2)!!} (u_{00} +  u_{01} + u_{10} + u_{11})\]
with
\[u_{00} = \frac{(n+6)!!}{(n+7)!!} \ , \ u_{01}=u_{10} = - 2 \frac{(n+4)!!}{(n+7)!!} \ , \ u_{11} = 4 \frac{(n+2)!!}{(n+7)!!}\]
Now, from Theorem 1.2 therein
\[I_{n-1} (2,2)  =\frac{(n-2)!!}{(n+2)!!}\]
which gives
\[ I_n\begin{pmatrix} 2&2\\2&2\end{pmatrix} = \frac{ (n^2+4n+7)}{(n+6) (n+4) (n+2)(n+1)n(n-1)}= \frac{1}{n^4} + O(\frac{1}{n^5})\]
The last one is managed with the same scheme:
\[ I_n\begin{pmatrix} 2&2\\0&2\end{pmatrix} = I_{n-1}(2,2) I_{n-1}(2) \Phi \begin{pmatrix} 2&2\\0&2\end{pmatrix} \]
and from Theorem 4.3 therein
\[\Phi \begin{pmatrix} 2&2\\0&2\end{pmatrix} = \frac{(n-1)!![(n+2)!!]^2}{(n-2)!! n!! (n+5)!!}\]
and since $I_{n-1}(2) = \mathbb E O_{11}^2= \frac{1}{n}$ we have got
\[ I\begin{pmatrix} 2&2\\0&2\end{pmatrix}  = \frac{(n-1)!!(n+2)!!}{n (n+2)!! n!! (n+5)!!}= \frac{(n+1)}{n^2 (n+4)(n+2)}= \frac{1}{n^3} + O(\frac{1}{n^4})\]
$\Box$
 
 \begin{lem} \label{lem1}
\begin{enumerate}
\item In the unitary case
 \begin{equation}
 \lim_{n \rightarrow \infty} n^2(n-1)^2 \mathbb E \left(V_{11}^2 V_{22}^2\right)  = 1
 \end{equation}
 and
 \begin{equation}
 \lim_{n \rightarrow \infty} n^2(n-1)^2 \mathbb E \left(V_{11}V_{12}V_{21}V_{22}\right)  = 0
 \end{equation}
\item In the orthogonal case
 \begin{equation}
 \lim_{n \rightarrow \infty} n^2(n-1)^2 \mathbb E \left(V_{11}^2 V_{22}^2\right)  = 4
 \end{equation}
 and
 \begin{equation}
 \lim_{n \rightarrow \infty} n^2(n-1)^2 \mathbb E \left(V_{11}V_{12}V_{21}V_{22}\right)  = 0
 \end{equation}
\end{enumerate}
 \end{lem}
 \proof
 We develop
 $$V_{11}^2 V_{22}^2 = (|U_{11}|^4 - \frac{2}{n}  |U_{11}|^2 + \frac{1}{n^2})(|U_{22}|^4 - \frac{2}{n}  |U_{22}|^2 + \frac{1}{n^2})$$ to obtain
 \begin{eqnarray}
\nonumber
 \lefteqn{
\E(V_{11}^2 V_{22}^2) = \E( |U_{11}|^4|U_{22}|^4) - \frac{4}{n} \E( |U_{11}|^4|U_{22}|^2) + \frac{2}{n^2} \E( |U_{11}|^4) } \\
\label{ueto}
&& + \frac{4}{n^2} \E( |U_{11}|^2|U_{22}|^2) -  \frac{4}{n^3} \E( |U_{11}|^2) +  \frac{1}{n^4}\,,
\end{eqnarray}
in the unitary case and the same expression with $U_{ij}$  replaced  by $O_{ij}$ in the orthogonal case. It remains to make the substitutions from Proposition \ref{propu}.
The proof of the second limit is similar.
$\Box$

\section{Appendix : Proofs of Theorem \ref{Lin1} and Theorem \ref{Lin2}}
 We define $D_n= \sum_{ij} |U_{ij}|^2 R_i Y_j$ and we prove that $A_n$ and $D_n$ have the same limit in law. We can write
$$A_n : =\sum_{j=1}^n \Lambda_j(n) C_j:= S_{n-1} + \Lambda_n(n) C_n , \quad  D_n  := \sum_{j=1}^n \Lambda_j(n) Y_j $$
where $\Lambda_j(n) = \sum_{i=1}^n  |U_{ij}|^2 R_i$ are independent of $C_j$ and $Y_j$.
Let $F$ be a smooth function with a bounded third derivative. The first step consists in replacing $C_n$ by $Y_n$ and to compare 
$\E(F( S_{n-1} +  \Lambda_n(n) C_n))$ and $\E(F( S_{n-1} +  \Lambda_n(n) Y_n))$.
Using a Taylor expansion,
\begin{eqnarray*}
\lefteqn{
\E(F( S_{n-1} +  \Lambda_n(n) C_n)) = \E(F( S_{n-1} ) )+  \E(F'( S_{n-1})) \E(\Lambda_n(n))\E(C_n)} \\
&& +  \frac{1}{2} \E(F''( S_{n-1})) \E(\Lambda_n(n)^2)\E(C_n^2) + O( \E(\Lambda_n(n)^3))\E(C_n^3)
\end{eqnarray*}
and a similar expression for $\E(F( S_{n-1} +  \Lambda_n(n) Y_n))$. All the terms in the two expressions, but the last, are equal. We thus need to estimate $\E(\Lambda_n(n)^3)$
By centering,
\begin{equation} \label{A2}
 \E(\Lambda_n(n)^3)= \sum_{i=1}^n \E(|U_{in}|^6 ) \E(R_i^3) = O(\frac{1}{n^2}) 
 \end{equation} 
from \eqref{A1} or \eqref{A1o}.
We repeat the operation of replacement of $C_j$ by $Y_j$ from $j=n-1$ to $1$  and by a summation, we obtain
$$ \E(F(A_n)) - \E(F(D_n)) = O(\frac{1}{n}).$$
In the same way $D_n$ and $B_n$ have the same limit in law, by exchanging the role of $i$ and $j$.  

We can extend this proof for the finite dimensional distributions. 
The proof is the analogue  as above, using a Taylor expansion of the two expressions, involving a smooth  function $F$ of $k$ variables,
$$ F(S_{n-1}(s_1,t_1) + \Lambda_n(s_1) (1_{(C_n \leq t_1)} -t_1), \ldots, S_{n-1}(s_k,t_k) + \Lambda_n(s_k) (1_{(C_n \leq t_k)} -t_k))$$
and
$$ F(S_{n-1}(s_1,t_1) + \Lambda_n(s_1) \gamma_n(t_1), \ldots, S_{n-1}(s_k,t_k) + \Lambda_n(s_k) \gamma_n(t_k)).$$
We then use that  the process $(1_{(C_i \leq t)} -t)$ has the same covariance as the Brownian bridge and the estimate \eqref{A2}.$\Box$

\renewcommand{\refname}{References}

\end{document}